\documentclass[nosumlimits]{amsart}

\usepackage{amscd,amssymb}

\theoremstyle{plain}
\newtheorem{prop}[subsection]{Proposition}
\newtheorem{thm}[subsection]{Theorem}
\newtheorem{lem}[subsection]{Lemma}
\newtheorem{cor}[subsection]{Corollary}

\theoremstyle{definition}
\newtheorem{defn}[subsection]{Definition}
\newtheorem{rem}[subsection]{Remark}

\numberwithin{equation}{section}

\newcommand{\DS}{\displaystyle }
\newcommand{\A}{{\mathcal A}}

\newcommand{\EE}{{\mathcal E}}
\newcommand{\Ii}{{\mathcal I}}
\newcommand{\JJ}{{\mathcal J}}
\newcommand{\LL}{{\mathcal L}}
\newcommand{\RR}{{\mathcal R}}
\newcommand{\Z}{{\mathbb Z}}
\newcommand{\C}{{\mathbb C}}
\newcommand{\T}{{({\mathbb C}^*)^N}}

\newcommand{\II}{{\mathbb I}}
\newcommand{\bC}{{\mathbf C}}
\newcommand{\hC}{{\widehat{C}}}
\newcommand{\hp}{{\widehat{\partial}}}

\newcommand{\bD}{{\mathbf D}}
\newcommand{\bM}{{\mathbf M}}
\newcommand{\bone}{{\mathbf 1}}
\newcommand{\bz}{{\mathbf z}}
\newcommand{\bt}{{\mathbf t}}
\newcommand{\bl}{{\boldsymbol{\lambda}}}
\newcommand{\bDD}{{\boldsymbol{\Delta}}}
\newcommand{\bX}{{\boldsymbol{\Xi}}}
\newcommand{\bp}{{\boldsymbol{\partial}}}
\newcommand{\bd}{{\boldsymbol{\delta}}}
\newcommand{\bc}{{\boldsymbol{\gamma}}}

\newcommand{\bzeta}{{\boldsymbol{\zeta}}}
\newcommand{\br}{{\boldsymbol{\rho}}}
\newcommand{\bJ}{{\boldsymbol{\mathcal J}}}
\newcommand{\D}{{\Delta }}
\newcommand{\p}{\partial }
\newcommand{\ii}{\text{i}}
\newcommand{\tr}{\tilde{\rho }}

\renewcommand{\b}{{\beta }}
\renewcommand{\c}{{\gamma }}

\renewcommand{\l}{{\lambda}}
\renewcommand{\L}{{\Lambda }}
\renewcommand{\ll}{{\ell }}
\renewcommand{\*}{{\bullet}}

\DeclareMathOperator{\rank}{rank}
\DeclareMathOperator{\id}{id}
\DeclareMathOperator{\Aut}{Aut}
\DeclareMathOperator{\End}{End}
\DeclareMathOperator{\Mat}{Mat}
\DeclareMathOperator{\homR}{Hom}
\DeclareMathOperator{\Hom}{{\mathcal H}\!\textit{\large om}}
\DeclareMathOperator{\TT}{T}

\begin{document}

\rightline{{\footnotesize In: Singularities and Arrangements, Sapporo-Tokyo 1998,}}
\rightline{{\footnotesize Advanced Studies in Pure Mathematics 
(to appear).}}
\bigskip
\bigskip

\title[Cohomology of Discriminantal Arrangements]
{On the Cohomology of Discriminantal Arrangements
and Orlik-Solomon Algebras}

\author[D.~Cohen]
{Daniel C.~Cohen}
\address{Department of Mathematics,
Louisiana State University,
Baton Rouge, LA 70803}
\email{cohen@math.lsu.edu}
\urladdr{http://math.lsu.edu/\~{}cohen}

\dedicatory{For Peter Orlik on the occasion of his sixtieth birthday.}

\thanks{Partially supported by
grant LEQSF(1996-99)-RD-A-04 from the Louisiana Board of Regents,}
\thanks{and by an NSF/EPSCoR Travel Grant for Emerging Faculty.}

\subjclass{Primary 52B30;  Secondary 20F36}

\keywords{discriminantal arrangement, local system, Orlik-Solomon 
algebra}

\begin{abstract} 
We relate the cohomology of the Orlik-Solomon algebra of a 
discriminantal arrangement to the local system cohomology of the 
complement.  The Orlik-Solomon algebra of such an arrangement (viewed 
as a complex) is shown to be a linear approximation of a complex 
arising from the fundamental group of the complement, the cohomology 
of which is isomorphic to that of the complement with coefficients in 
an arbitrary complex rank one local system.  We also establish the 
relationship between the cohomology support loci of the complement of 
a discriminantal arrangement and the resonant varieties of its 
Orlik-Solomon algebra.
\end{abstract}

\maketitle

\section*{Introduction}
\label{sec:intro}

Let $\A$ be an arrangement of $N$ complex hyperplanes, and let $M(\A)$ 
be its complement.  For each hyperplane $H$ of $\A$, let $f_H$ be a 
linear polynomial with kernel $H$, and let $\l_H$ be a complex number.  
Each point $\bl=(\dots,\l_H,\dots)\in\C^N$ determines an integrable 
connection $\nabla=d+\Omega_{\bl}$ on the trivial line bundle over 
$M(\A)$, where $\Omega_{\bl}=\sum_{H\in\A}\l_Hd\log f_H$, and an 
associated complex rank one local system $\LL$ on $M(\A)$.  
Alternatively, if $\bt\in\T$ is the point in the complex torus 
corresponding to $\bl$, then the local system $\LL$ is induced by the 
representation of the fundamental group of $M(\A)$ which sends any 
meridian about $H\in\A$ to $t_H=\exp(-2\pi\ii\l_H)$.

Due largely to its various applications, the cohomology of $M(\A)$ 
with coefficients in $\LL$ has been the subject of considerable recent 
interest.  These applications include representations of braid groups, 
generalized hypergeometric functions, and the Knizhnik-Zamolodchikov 
equations from conformal field theory.  See, for instance, the works 
of Aomoto, Kita, Kohno, Schechtman, and 
Varchenko~\cite{Ao,AK,Ko,SV,Va}, 
and see Orlik and Terao~\cite{OT} as a general reference for 
arrangements.  Of particular interest in these applications are the 
discriminantal arrangements of \cite{SV}, the complements of which 
may be realized as configuration spaces of 
ordered points in $\C$ punctured finitely many times.
(Note that our use of the term 
``discriminantal'' differs from that of \cite{OT}.)

The local system cohomology $H^*(M(\A);\LL)$ may be studied from a 
number of points of view.  For instance, if $\A$ is real, that is, 
defined by real equations, the complement $M(\A)$ is homotopy 
equivalent to the Salvetti complex $X$ of $\A$, see~\cite{Sal}.  In 
this instance, the complex $X$ may be used in the study of local 
systems on $M(\A)$.  This approach is developed by Varchenko in 
\cite{Va}, to which we also refer for discussion of the applications 
mentioned above, and has been pursued by Denham and Hanlon~\cite{DH} 
in their study of the homology of the Milnor fiber of an arrangement.

If $\A$ is $K(\pi,1)$, that is, the complement $M(\A)$ is a 
$K(\pi,1)$-space, then local systems on $M(\A)$ may be studied from 
the point of view of cohomology of groups.  Any representation of the 
fundamental group $G$ of the complement of a $K(\pi,1)$ arrangement 
gives rise to a $G$-module $L$, and a local system of coefficients 
$\LL$ on $M(\A)$.  Since $M(\A)$ is a $K(\pi,1)$-space, we have 
$H_{*}(M(\A);\LL)=H_{*}(G;L)$ and $H^{*}(M(\A);\LL)=H^{*}(G;L)$, see 
for instance Brown \cite{Brown}.  The class of $K(\pi,1)$ arrangements 
includes the discriminantal arrangements noted above, as they are 
examples of fiber-type arrangements, well-known to be $K(\pi,1)$, see 
e.g.~Falk and Randell~\cite{FR1}.

For any arrangement $\A$, let $B(\A)$ denote the Brieskorn algebra of 
$\A$, generated by $1$ and the closed differential forms $d\log f_H$, 
$H\in\A$.  As is well-known, the algebra $B(\A)$ is isomorphic to 
$H^*(M(\A);\C)$, and to the Orlik-Solomon algebra $A(\A)$, so is 
determined by the lattice of $\A$, see \cite{Br,OS,OT}.  If $\LL$ is a 
local system on $M(\A)$ determined by ``weights'' $\bl$ which satisfy 
certain Aomoto non-resonance conditions, work of Esnault, Schechtman, 
and Viehweg~\cite{ESV}, extended by Schechtman, Terao, and 
Varchenko~\cite{STV}, shows that $H^*(M(\A);\LL)$ is isomorphic to the 
cohomology of the complex $(B(\A),\Omega_{\bl}\wedge)$.  Thus for 
non-resonant weights, the local system cohomology may be computed by 
combinatorial means, using the Orlik-Solomon algebra equipped with 
differential $\mu(\bl)$, given by left-multiplication by 
$\omega_{\bl}$, the image of $\Omega_{\bl}$ under the isomorphism 
$B(\A)\to A(\A)$.

For arbitrary (resonant) weights, one has
\begin{equation*} 
\dim H^k(A(\A),\mu(\bl)) \le \dim H^k(M(\A);\LL) \le \dim 
H^k(M(\A);\C)
\end{equation*}
for each $k$.  See Libgober and Yuzvinsky~\cite{LY} for the first of 
these inequalities.  The second is obtained using stratified Morse 
theory in~\cite{morse}, and resolves a question raised by Aomoto and 
Kita in~\cite{AK}.  For resonant weights, the precise relation between 
$H^*(A(\A),\mu(\bl))$ and $H^*(M(\A);\LL)$ is not known.

However, recent results suggest, at least for small $k$, that 
$H^k(A(\A),\mu(\bl))$ may be viewed as a ``linear approximation'' of 
$H^k(M(\A);\LL)$.  The resonant varieties, 
$\RR_k^{m}(A(\A))=\{\l\in\C^N \mid \dim H^k(A(\A),\mu(\bl)) \ge m\}$, 
of the Orlik-Solomon algebra were introduced by Falk in \cite{Fa}.  
For $k=1$ and any arrangement $\A$, it is known that 
$\RR_{1}^{m}(A(\A))$ coincides with the tangent cone of the 
cohomology support locus of the complement, 
$\Sigma^{1}_m(M(\A))=\{\bt\in\T\mid \dim H^1(M(\A);\LL) \ge m\}$, 
at the point $(1,\dots,1)$, see \cite{CScv, L3, LY}.  For certain 
arrangements, we present further ``evidence'' in support of this 
philosophy~here.

If $\A$ is a fiber-type arrangement, the fundamental group $G$ of the 
complement $M(\A)$ may be realized as an iterated semidirect product 
of free groups, and $M(\A)$ is a $K(G,1)$-space, see \cite{FR1,OT}.  
For any such group, we construct a finite, free $\Z{G}$-resolution, 
$C_\*(G)$, of $\Z$ in \cite{CScc}.  This resolution may be used to 
compute the homology and cohomology of $G$ with coefficients in any 
$G$-module $L$, or equivalently, that of $M(\A)$ with coefficients in 
any local system $\LL$.  We have $H_*(M(\A);\LL) = 
H_*(C_\*(G)\otimes_G L)$ and $H^*(M(\A);\LL) = 
H^*(\Hom_G(C_\*(G),L))$, see \cite{Brown}.

Briefly, for a fiber-type arrangement $\A$, the relationship between 
the cohomology theories $H^*(A(\A),\mu(\bl))$ and $H^*(M(\A);\LL)$ is 
given by the following assertion.  {\em For any $\bl$, the complex 
$(A(\A),\mu(\bl))$ is a linear approximation of the complex 
$\Hom_G(C_\*(G),L)$.} We prove a variant of this statement in the case 
where $\A$ is a discriminantal arrangement here.  We also establish 
the relationship between the resonant varieties $\RR_k^{m}(A(\A))$ and 
cohomology support loci $\Sigma_m^{k}(M(\A))$ of these arrangements, 
analogous to that mentioned above in the case $k=1$.

The paper is organized as follows.  The Orlik-Solomon algebra of a 
discriminantal arrangement admits a simple description, which 
fascilitates analysis of the differential of the complex 
$(A(\A),\mu(\bl))$.  We carry out this analysis, which is elementary 
albeit delicate, in section~\ref{sec:H*OS}, and obtain an explicit 
(inductive) description of the differential $\mu(\bl)$.  In section 
\ref{sec:H*G}, we recall the construction of the resolution 
$C_{\*}(G)$ from \cite{CScc} in the instance where $G$ is the 
fundamental group of the complement of a discriminantal arrangement, 
and exhibit a complex $(\bC^\*,\bd^\*(\bt))$ which computes the cohomology 
$H^*(M(\A),\LL)$ for an arbitrary rank one local system.  
We then study in section \ref{sec:dercx} a linear approximation of 
$(\bC^\*,\bd^\*(\bt))$, and relate it, for arbitrary $\bl$, to the complex 
$(A(\A),\mu(\bl))$.  We conclude by realizing the resonant varieties 
of the Orlik-Solomon algebra of a discriminantal arrangement as the 
tangent cones at the identity of the cohomology support loci of the 
complement in section \ref{sec:loci}.

\section{Cohomology of the Orlik-Solomon Algebra}
\label{sec:H*OS}
Let $M_{n}=\{(x_{1},\dots,x_{n})\in\C^{n}\mid x_{i}\neq x_{j} \text{ 
if } i\neq j\}$ be the configuration space of $n$ ordered points in 
$\C$.  Note that $M_{n}$ may be realized as the complement of the 
braid arrangement $\A_{n}=\{x_{i}=x_{j}, 1\le i<j\le n\}$ in $\C^{n}$.  
Classical work of Fadell and Neuwirth~\cite{FN} shows the projection 
$\C^n \to \C^\ell$ defined by forgetting the last $n-\ell$ coordinates 
gives rise to a bundle map $p:M_n \to M_\ell$.  From this it follows 
that $M_n$ is a $K(P_n,1)$-space, where $P_n = \pi_1(M_n)$ is the pure 
braid group on $n$ strands.

The typical fiber of the bundle of configuration spaces $p:M_n\to 
M_\ell$ may be realized as the complement of an arrangement in 
$\C^{n-\ell}$, a discriminantal arrangement in the sense of Schechtman 
and Varchenko, see \cite{SV,Va}.  The fiber over $\bz=(z_{1},\dots, 
z_{\ll})\in M_{\ell}$ may be realized as the complement, 
$M_{n,\ll}=M(\A_{n,\ll})$, of the arrangement $\A_{n,\ell}$ consisting 
of the $N=\binom{n}{2}-\binom{\ell}{2}$ hyperplanes
\begin{equation*}
H_{i,j}=
\begin{cases}
\ker(x_j-x_i)&\ell+1 \le i < j \le n,\\
\ker(x_j-z_{i})&1\le i \le \ell,\ \ell+1\le j \le n,
\end{cases}
\end{equation*}
in $\C^{n-\ell}$ (with coordinates $x_{\ell+1},\dots,x_n$).  Note that 
$M_{n,\ll}$ is the configuration space of $n-\ell$ ordered points in 
$\C\setminus\{z_1,\dots,z_\ll\}$, and that the topology of $M_{n,\ll}$ 
is independent of $\bz$, see \cite{FN,Bi,Ha}.  We first record some 
known results on the cohomology of $M_{n,\ell}$.

\subsection{The Orlik-Solomon Algebra}
\label{subsec:Hbraid}
The fundamental group of $M_{n,\ll}$ may be realized as 
$P_{n,\ell}=\pi_1(M_{n,\ell})=\ker(P_n\to P_\ell)$, the kernel of the 
homorphism from $P_{n}$ to $P_{\ll}$ defined by forgetting the last 
$n-\ll$ strands.  From the homotopy exact sequence of the bundle 
$p:M_n\to M_\ell$, we see that $M_{n,\ell}$ is a 
$K(P_{n,\ell},1)$-space.  The cohomology of this space, and hence of 
this group, may be described as follows.

Let $\EE=\bigoplus_{q=0}^N \EE^q$ be the 
graded exterior algebra over $\C$, generated by 
$e_{i,j}$, $\ll+1\le j\le n$, $1\le i < j$.
Let $\Ii$ be the ideal 
in $\EE$ generated, for $1\le i<j<k\le n$,~by
\[ 
e_{i,j}\wedge e_{i,k}-e_{i,j}\wedge e_{j,k}+e_{i,k}\wedge e_{j,k} 
\text{ if $j\ge\ll+1$,\quad and }\quad 
e_{i,k}\wedge e_{j,k} \text{ if $j\le\ll$.}
\] 
Note that $\EE^q\subset\Ii$ for $q>n-\ll$.  
The Orlik-Solomon algebra of the discriminantal arrangement 
$\A_{n,\ll}$ is the quotient $A=\EE/\Ii$.

\begin{thm}\label{thm:H^*Pn}
The cohomology algebra $H^*(M_{n,\ell};\C)=H^*(P_{n,\ell};\C)$ is 
isomorphic to the Orlik-Solomon algebra $A=A(\A_{n,\ll})$.
\end{thm}

The grading on $\EE$ induces a grading 
$A=\bigoplus_{q=0}^{n-\ll}A^{q}$ on the Orlik-Solomon algebra 
$A=A(\A_{n,\ll})$.  Let $a_{i,j}$ denote the image of $e_{i,j}$ in 
$A$, and note that these elements form a basis for $A^{1}$ and 
generate $A$.  From the description of the ideal $\Ii$ above, it is 
clear that all relations among these generators are consequences of 
the following:
\begin{equation} \label{eqn:APnrels}
a_{i,k} \wedge a_{j,k}= 
\begin{cases}
a_{i,j} \wedge (a_{j,k} - a_{i,k}) &\text{if $j \ge \ell+1$,}\\
0 & \text{if $j \le \ell$,}
\end{cases}
\end{equation}
for $1 \le i<j<k \le n$.

This observation leads to a natural choice of basis for the algebra 
$A$.  For $m\le n$, write $[m,n]=\{m,m+1,\dots,n\}$.  If 
$I=\{i_1,\dots,i_q\}$ and $J=\{j_1,\dots,j_q\}$ satisfy the conditions 
$J\subseteq [\ll+1,n]$ and $1\le i_p < j_p$ for each $p$, let 
$a_{I,J}=a_{i_1,j_1}\wedge\dots\wedge a_{i_q,j_q}$.  If $|J|=0$, set 
$a_{I,J}=1$.

\begin{prop} \label{prop:Abasis}
For each $q$, $0\le q\le n-\ell$, the forms $a_{I,J}$ with $|J|=q$ and 
$I$ as above form a basis for the summand $A^q$ of the Orlik-Solomon 
algebra $A$ of the discriminantal arrangement $\A_{n,\ll}$.  
Furthermore, the summand $A^q$ decomposes as a direct sum, 
$A^q=\bigoplus_{|J|=q} A_J$, where $A_J = \bigoplus_{I} \C\, a_{I,J}$.
\end{prop}

\begin{rem} These results are well-known.  For instance, if $\A=\A_n$ 
is the braid arrangement, Theorem \ref{thm:H^*Pn} follows from results 
of Arnol'd \cite{Arn} and Cohen \cite{fred}, which show that 
$H^{*}(M_{n};\C)$ is generated by the forms 
$a_{i,j}=d\log(x_{j}-x_{i})$, with relations \eqref{eqn:APnrels} (with 
$\ll=1$).  For any discriminantal arrangement $\A_{n,\ll}$, Theorem 
\ref{thm:H^*Pn} is a consequence of results of Brieskorn and 
Orlik-Solomon, see \cite{Br,OS,OT}.

As mentioned in the introduction, the discriminantal arrangements 
$\A_{n,\ll}$ are examples of (affine) fiber-type or supersolvable 
arrangements.  The structure of the Orlik-Solomon algebra of any such 
arrangement $\A$ was determined by Terao \cite{hiro}.  The basis for 
the algebra $A(\A_{n,\ll})$ exhibited in Proposition \ref{prop:Abasis} 
above is the \textbf{nbc}-basis (with respect to a natural ordering of 
the hyperplanes of $\A_{n,\ll}$), see \cite{OT}.  The Orlik-Solomon 
algebra of any supersolvable arrangement admits an analogous basis, 
see Bj\"orner-Ziegler \cite{BZ}, and see Falk-Terao \cite{FT} for 
affine supersolvable arrangements.
\end{rem}

\subsection{The Orlik-Solomon Algebra as a Complex}
\label{subsec:wedge}
Recall that $N=\binom{n}{2}-\binom{\ell}{2}$, and consider $\C^N$ with 
coordinates $\l_{i,j}$, $\ell+1\le j \le n$, $1\le i <j$.  Each point 
$\bl \in \C^N$ gives rise to an element $\omega=\omega_\bl = \sum 
\l_{i,j}\cdot a_{i,j}$ of $A^1$.  Left-multiplication by $\omega$ 
induces a map $\mu^q(\bl):A^{q}\to A^{q+1}$, defined by 
$\mu^q(\bl)(\eta)=\omega\wedge\eta$.  Clearly, $\mu^{q+1}(\bl) \circ 
\mu^q(\bl)=0$, so $(A^\*,\mu^\*(\bl))$ is a complex.

We shall obtain an inductive formula for the boundary maps of the 
complex $(A^{\*},\mu^\*(\bl))$.  The projection $\C^{n-\ll}\to\C$ onto 
the first coordinate gives rise to a bundle of configuration spaces, 
$M_{n,\ell} \to M_{\ll+1,\ell}$, with fiber $M_{n,\ell+1}$, see 
\cite{FN,Bi,Ha}.  The inclusion of the fiber 
$M_{n,\ell+1}\hookrightarrow M_{n,\ell}$ induces a map on cohomology 
which is clearly surjective.  This yields a surjection 
$\pi:A(\A_{n,\ell}) \to A(\A_{n,\ell+1})$.

Write $A=A(\A_{n,\ell})$ and $\widehat A=A(\A_{n,\ell+1})$, and denote 
the generators of both $A$ and $\widehat A$ by $a_{i,j}$.  In terms of 
these generators, the map $\pi$ is given by $\pi(a_{i,\ll+1})=0$, and 
$\pi(a_{i,j})=a_{i,j}$ otherwise.  Let $\widehat\omega \in \widehat A$ 
denote the image of $\omega\in A^1$ under $\pi$.  If we write $\omega= 
\sum_{k=\ell+1}^n \omega_k$, where $\omega_k = \sum_{i=1}^{k-1} 
\l_{i,k} \cdot a_{i,k}$, then $\widehat\omega = \sum_{k=\ell+2}^n 
\omega_k$.  As above, left-multiplication by $\widehat\omega$ induces 
a map $\widehat\mu^q(\bl):\widehat A^q \to \widehat A^{q+1}$, and 
$(\widehat A^\*,\widehat\mu^\*(\bl))$ is a complex.  The following is 
straightforward.

\begin{lem} The map $\pi:(A^\*,\mu^\*(\bl)) \to (\widehat 
A^\*,\widehat\mu^\*(\bl))$ is a surjective chain map.
\end{lem}

Let $(B^\*,\mu_B^\*(\bl))$ denote the kernel of the chain map $\pi$.  
The terms are of the form $B^q = \bigoplus A^q_K$, where $\ll+1\in K$ 
and $|K|=q$.  In particular, $B^0=0$.  We now identify the 
differential $\mu_B^\*(\bl)$.  If $k<m\le n$ and $J\subseteq [m,n]$, 
let $\{k,J\}$ denote the (ordered) subset $\{k\}\cup J$ of $[k,n]$.  
For a linear map $F$, write $[F]^k$ for the map $\oplus_1^k F$.

\begin{prop} \label{prop:directsum}
The complex $(B^\*,\mu_B^\*(\bl))$ decomposes as the direct sum of 
$\ell$ copies of the complex $\widehat A^\bullet$, shifted in 
dimension by one, with the sign of the boundary map reversed.  In 
other words, $(B^\*,\mu_B^\*(\bl)) \cong \bigl((\widehat 
A^{\*-1})^\ll,-[\widehat\mu^{\*-1}(\bl)]^\ll\bigr)$.
\end{prop}
\begin{proof}
For $1\le q\le n-\ll$, we have $B^q=\bigoplus A^q_{\{\ell+1,J\}}$, 
where the sum is over all $J\subseteq [\ll+2,n]$ with $|J|=q-1$.  Each 
summand may be written as $A^q_{\{\ell+1,J\}} = \bigoplus_{i=1}^{\ell} 
a_{i,\ell+1} \wedge A^{q-1}_J$.  Thus, $B^q = \bigoplus_{i=1}^{\ell} 
a_{i,\ell+1}\wedge\widehat A^{q-1}$ is isomorphic to the direct sum of 
$\ell$ copies of $\widehat A^{q-1}$ via the map $B^{q}\to [\widehat 
A^{q-1}]^{\ll}$, $a_{i,\ll+1}\wedge a_{I,J} \mapsto 
(0,\dots,a_{I,J},\dots,0)$.

Now consider the boundary map $\mu_B^q(\bl):B^q\to B^{q+1}$ of the 
complex $B^\*$, induced by left-multiplication by 
$\omega=\sum_{k=\ell+1}^n\omega_k$.  Let $\eta=a_{i,\ell+1}\wedge 
a_{I,J}$ be a generator for $B^q$.  Since $a_{i,k}\wedge a_{j,k}=0$ 
for all $i,j<k$, we have $\omega_{\ell+1}\wedge \eta=0$.  Thus, 
\[\mu_B^q(\bl)(\eta)=\omega\wedge\eta= 
(\omega-\omega_{\ell+1})\wedge\eta = 
-a_{i,\ll+1}\wedge(\omega-\omega_{\ell+1})\wedge a_{I,J}.
\]
Write $(\omega-\omega_{\ell+1})\wedge a_{I,J} = \xi_1+\xi_2$ in terms 
of the basis for $A$ specified in Proposition~\ref{prop:Abasis}, where 
$\xi_1 \in \bigoplus_{\ell+1\in K} A^q_K$ and $\xi_2\in 
\bigoplus_{\ell+1\notin K} A^q_K$.  Then we have $ \omega\wedge\eta = 
-a_{i,\ell+1}\wedge(\xi_1+\xi_2)=-a_{i,\ell+1}\wedge\xi_2$.  Checking 
that $\widehat\omega\wedge a_{I,J} = \xi_2$ in $\widehat A$, we have 
$\mu_B^q(\bl)(a_{i,\ell+1}\wedge a_{I,J})=-a_{i,\ell+1}\wedge 
\widehat\mu^{q-1}(\bl)( a_{I,J})$.  Thus, with the change in sign, the 
boundary map $\mu_B^\*(\bl)$ respects the direct sum decomposition 
$B^\* \cong (\widehat A^{\*-1})^{\ll}$.
\end{proof}

\subsection{Boundary Maps}
\label{subsec:OSmaps}
We now study the differential of the complex $(A^{\*},\mu^{\*}(\bl))$.  
The direct sum decompositions of the terms of the complexes $A^{\*}$, 
$\widehat A^{\*}$, and $B^{\*}$ exhibited above yield 
\[
A^{q}=\bigoplus_{|J|=q} A^{q}_{J} =\Bigl(\bigoplus_{\ell+1\in J} 
A^{q}_{J}\Bigr) \oplus \Bigl(\bigoplus_{\ell+1\notin J} 
A^{q}_{J}\Bigr) =B^{q} \oplus \widehat A^{q}.
\]

Let $\pi_{B}:A^{q} \to B^{q}$ denote the natural projection.  With 
respect to the direct sum decomposition of the terms 
$A^{q}=B^{q}\oplus \widehat A^{q}$, the boundary map $\mu^{\*}(\bl)$ 
of the complex $A^{\*}$ is given by $\mu^{q}(\bl)(v_{1},v_{2}) = 
(\mu^{q}_{B}(\bl)(v_{1}) + \Psi^q(\bl)(v_{2}),\widehat 
\mu^{q}(\bl)(v_{2}))$, where $\Psi^{q}(\bl) = \pi_{B}\circ 
\mu^{q}(\bl):\widehat A^{q}\to B^{q+1}$.  In matrix form, we have
\begin{equation} \label{eq:osmapmatrix}
\mu^{q}(\bl) = \pmatrix \mu^{q}_{B}(\bl)&0\\
\Psi^q(\bl)&\widehat \mu^{q}(\bl)\endpmatrix.
\end{equation}

Since $\widehat A^{\*}$ is the complex associated to the 
discriminantal arrangement $\A_{n,\ell+1}$ in $\C^{n-\ll-1}$ and 
$B^{\*} \cong (\widehat A^{\*-1})^{\ell}$ decomposes as a direct sum 
by Proposition~\ref{prop:directsum}, we inductively concentrate our 
attention on the maps $\Psi^q(\bl)$.  Fix $J\subseteq [\ll+2,n]$, and 
denote the restriction of $\Psi^q(\bl)$ to the summand $A^q_J$ of 
$\widehat A^q$ by $\Psi^q_J(\bl)$.  For $\eta\in A^q_J$, since 
$\pi_B(\omega_k\wedge\eta)=0$ if $k\notin \{\ll+1,J\}$, we have 
$\Psi^q_J(\bl)(\eta)=\omega_{\ll+1}\wedge\eta+\sum_{j\in 
J}\pi_B(\omega_j \wedge\eta)$.  Thus, $\Psi^q_J(\bl):A^q_J\to 
A^{q+1}_{\{\ll+1,J\}}= \bigoplus_{m=1}^\ll a_{m,\ll+1}\wedge A^q_J$.

For $1\le m\le\ll$, let $\pi_{m,\ell+1}:A_{\{\ell+1,J\}}^{q+1} \to 
a_{m,\ell+1} \wedge A_J^q$ denote the natural projection.  Then (the 
matrix of) $\Psi^q_J(\bl): A^q_J\to\bigl(A^q_J\bigr)^\ll$ may be 
expressed as
\begin{equation} \label{eq:psiblocks}
\Psi^q_J(\bl)=\begin{pmatrix} 
\pi_{1,\ell+1}\circ\Psi^q_J(\bl) & \cdots &
\pi_{m,\ell+1}\circ\Psi^q_J(\bl) & \cdots &
\pi_{\ll,\ell+1}\circ\Psi^q_J(\bl)
\end{pmatrix}, 
\end{equation}
and we focus our attention on one such block, that is, on the 
composition
\begin{equation} \label{eq:wedgew}
\pi_{m,\ell+1}\circ\Psi^q_J(\bl):A_J^q \longrightarrow 
A_{\{1,J\}}^{q+1} \longrightarrow 
a_{m,\ell+1} \wedge A_J^q.
\end{equation}

Write $J=\{j_1,\dots,j_q\}$ and for $1\le p \le q$, let $J_p = 
\{j_1,\dots,j_p\}$ and $J^{p} = J \setminus J_p$.  If $p=0$, set 
$J_0=\emptyset$ and $J^{0}=J$.  Then for $a_{I,J}\in A^{q}_{J}$, it is 
readily checked that $\pi_{m,\ell+1}\circ\Psi^q_J(\bl)(a_{I,J})= 
\pi_{m,\ell+1}\circ\pi_B(\omega\wedge a_{I,J}) $ is given by
\begin{equation} \label{eq:reclem}
\pi_{m,\ell+1}\circ\Psi^q_J(\bl)(a_{I,J}) =
\sum_{p=0}^q 
\pi_{m,\ell+1}\circ\pi_{B}(\omega_{j_p}\wedge a_{I_p^{},J_p^{}})
\wedge a_{I^{p}_{},J^{p}_{}},
\end{equation}
where $j_0=\ll+1$.  In light of this, we restrict our attention to 
$\pi_{m,\ll+1}\circ\pi_B(\omega_{j_q}\wedge a_{I,J})$.  We describe 
this term using the following notion.

\begin{defn} \label{def:admissible} 
Fix $J=\{j_1,\dots,j_q\}\subseteq[\ll+2,n]$ and $m\le\ll$.  If 
$I=\{i_1,\dots,i_q\}$ and $1\le i_p < j_p$ for each $p$, a set 
$K=\{k_{s_1},\dots,k_{s_t},k_{s_{t+1}}\}$ is called {\em 
$I$-admissible} if
\begin{enumerate}
\item $\{i_{s_1},\dots,i_{s_t}\} \subseteq I \setminus \{i_q\}$ 
and $i_{s_{t+1}} = i_q$;
\item $\{k_{s_1},i_{s_1}\} = \{m,\ll+1\}$;  and
\item $\{k_{s_p},i_{s_p}\} = \{k_{s_{p-1}},j_{s_{p-1}}\}$ for 
$p=2,\dots,t+1$.
\end{enumerate}
Note that the last condition is vacuous if $K$ is of cardinality one.  
Note also that $1\le k_{s_p} < j_{s_p}$ and $k_{s_p}\neq i_{s_p}$ for 
each $p$.
\end{defn}

\begin{lem} \label{lem:summands} We have
\[
\pi_{m,\ll+1}\circ\pi_B(\omega_{j_q}\wedge a_{I,J}) = 
\sum_K \l_{k_q,j_q} a_{m,\ll+1}\wedge b_{j_1}\wedge
\dots \wedge b_{j_q},
\]
where the sum is over all $I$-admissible sets $K=\{k_{s_{1}},\dots, 
k_{s_{t}},k_{s_{t+1}}=k_{q}\}$, and 
\[
b_{j_p} = 
\begin{cases}
a_{i_p,j_p} - a_{k_p,j_p}&\text{if $p \in \{s_{1},\dots,s_{t},q\}$,}\\
a_{i_p,j_p}&\text{if $p \notin \{s_{1},\dots,s_{t},q\}$.}
\end{cases}
\]
\end{lem}
\begin{proof} Let $a_{i,j}$ and $a_{k,j}$ be elements of 
$A^1_{\{j\}}$.  Write $r=\min\{i,k\}$ and $s=\max\{i,k\}$.  From 
\eqref{eqn:APnrels}, we have either $a_{i,j} \wedge a_{k,j} = 
a_{r,s}\wedge (a_{k,j} - a_{i,j})$ if $s\ge\ll+1$, or $a_{i,j} \wedge 
a_{k,j}=0$ if $s\le\ll$.  It follows from these considerations, and a 
routine exercise to check the sign, that summands $\l_{k_q,j_q} 
a_{m,\ll+1}\wedge b_{j_1}\wedge \dots \wedge b_{j_q}$ of 
$\pi_{m,\ll+1}\circ\pi_B(\omega_{j_q}\wedge a_{I,J})$ arise only from 
$I$-admissible sets $K$.
\end{proof}

Now write $\pi_{m,\ll+1}\circ\pi_B(\omega_{j_q}\wedge a_{I,J})= \sum_R 
\l^J_{R,I} a_{m,\ll+1}\wedge a_{R,J}$, where the sum is over all 
$R=\{r_1,\dots,r_q\}$, $1\le r_p < j_p$, $1\le p \le q$, and 
$\l^J_{R,I}\in \C$.

\begin{prop} \label{prop:wedgewjq}
The coefficient $\l^J_{R,I}$ of $a_{m,\ll+1}\wedge a_{R,J}$ in
$\pi_{m,\ll+1}\circ\pi_B(\omega_{j_q}\wedge a_{I,J})$ is given by
\[
\l^J_{R,I} = (-1)^{|R\setminus R\cap I|} \sum_K \l_{k_q,j_q}
\]
where the sum is over all $I$-admissible sets $K$ such that 
$R \setminus R\cap I \subseteq K$.
\end{prop}
\begin{proof}
Let $K=\{k_{s_1},\dots,k_{s_t},k_q\}$ be an $I$-admissible set.  
Associated with $K$, we have the term $\l_{k_q,j_q} a_{m,\ll+1}\wedge 
b_{j_1}\wedge \dots \wedge b_{j_q}$ of 
$\pi_{m,\ll+1}\circ\pi_B(a_{I,J}\wedge\omega_{j_q})$ from 
Lemma~\ref{lem:summands}.  If $R\setminus R\cap I \not\subseteq K$, it 
is readily checked that this term contributes nothing to the 
coefficient $\l^J_{R,I}$ of $a_{m,\ll+1}\wedge a_{R,J}$.  On the other 
hand, if $R\setminus R\cap I \subseteq K$, then the above term 
contributes the summand $(-1)^{|R\setminus R\cap I|} \l_{k_q,j_q}$ to 
the coefficient $\l^J_{R,I}$.
\end{proof}

We now obtain a complete description of the map $\pi_{m,\ll+1} \circ 
\Psi^q_J(\bl):A^q_J \to a_{m,\ll+1}\wedge A^q_J$ 
from~\eqref{eq:wedgew}.  Write 
$\pi_{m,\ll+1}\circ\Psi^q_J(\bl)(a_{I,J}) = \sum_R \L^J_{R,I} 
a_{m,\ll+1}\wedge a_{R,J}$, where, as above, the sum is over all 
$R=\{r_1,\dots,r_q\}$, $1\le r_p <j_p$, $1\le p \le q$, and 
$\L^J_{R,I}\in \C$.  Let $\epsilon_{R,I}=1$ if $R=I$, and 
$\epsilon_{R,I}=0$ otherwise.

\begin{thm} \label{thm:OSformula}
The coefficient $\L^J_{R,I}$ of $a_{m,\ll+1}\wedge a_{R,J}$ in
$\pi_{m,\ll+1} \circ \Psi^q_J(\bl)$ is given~by
\[
\L^J_{R,I} = (-1)^{|R\setminus R \cap I|}\Bigl( 
\epsilon_{R,I} \l_{m,\ll+1} +
\sum_{j\in J} \sum_K \l_{k,j}\Bigr),
\]
where, if $j=j_p$, the second sum is over all $I_p$-admissible sets 
$K=\{k_{s_1},\dots,k_{s_t},k\}$ for which $R \setminus R\cap I 
\subseteq K$.
\end{thm}
\begin{proof} From~\eqref{eq:reclem}, we have 
\[  
\pi_{m,\ll+1}\circ\Psi^q_J(a_{I,J}) =
\sum_{p=0}^q \pi_{m,\ll+1}\circ\pi_{B}(
\omega_{j_p}\wedge a_{I_p,J_p})
\wedge a_{I^{^p},J^{^p}},
\]  
and the summand corresponding to $p=0$ is simply $\l_{m,\ll+1} 
a_{m,\ll+1}\wedge a_{I,J}$.  For $p \ge 1$, write 
$\pi_{m,\ll+1}\circ\pi_{B}(\omega_{j_p}\wedge a_{I_p,J_p})= \sum_{R_p} 
\l^{J_p}_{R_p,I_p} a_{m,\ll+1}\wedge a_{R_p,J_p}$, where the sum is 
over all $R_p=\{r_1,\dots,r_p\}$.  For a fixed $R$, the coefficient of 
$a_{m,\ll+1} \wedge a_{R,I}$ in $\pi_{m,\ll+1}\circ\pi_B(\omega\wedge 
a_{I,J})$ may then be expressed as 
\[
\L^J_{R,I} = \epsilon_{R,I} \l_{m,\ll+1} + \sum_{p=1}^q 
\l^{J_p}_{R_p,I_p},
\]
where $R=R_p \cup I^p$.  Note that we have $R\setminus R\cap I=
R_{p}\setminus R_p \cap I_p$ for such $R$.

By Proposition~\ref{prop:wedgewjq}, we have $\l^{J_p}_{R_p,I_p} = 
(-1)^{|R_{p}\setminus R_p\cap I_p|} \sum_K \l_{k_p,j_p}$, where the 
sum is over all $I_p$-admissible sets $K$ with $R_p \setminus R_p\cap 
I_p \subseteq K$.  Thus,
\[
\L^J_{R,I} = \epsilon_{R,I} \l_{m,\ll+1} + \sum_{p=1}^q 
(-1)^{|R_{p}\setminus R_p\cap I_p|} \sum_K \l_{k_p,j_p},
\]
and since $R=R_p\cup I^p$, we have $R\setminus R\cap I=
R_p\setminus R_p\cap I_p \subseteq K$.
\end{proof}

\begin{rem} In light of the decomposition of the boundary maps of the 
complex $(A^\*,\mu^\*(\bl))$ given by \eqref{eq:osmapmatrix} and 
\eqref{eq:psiblocks}, the above theorem, together with the ``initial 
conditions'' $\mu^0(\bl):A^0\to A^1$, $1\mapsto \sum_{k=\ll+1}^n 
\omega_k=\sum_{k=\ll+1}^n \sum_{i=1}^{k-1} \l_{i,k}a_{i,k}$, provides 
a complete description of the boundary maps $\mu^\*(\bl)$.
\end{rem}

\section{Resolutions and Local Systems}
\label{sec:H*G}
The fundamental group of the complement of a discriminantal 
arrangement, and more generally that of any fiber-type arrangement, 
may be realized as an iterated semidirect product of free groups.  For 
any such group $G$, in \cite{CScc} we construct a finite free $\Z{G}$ 
resolution $C_\*(G)$ of the integers.  We recall the construction of 
this resolution in notation consonant with that of the previous 
section.

Denote the standard generators of the pure braid group $P_n$ by 
$\c_{i,j}$, $1\le i<j\le n$, and for each $j$, let $G_j$ be the free 
group on the $j-1$ generators $\c_{1,j},\dots,\c_{j-1,j}$.  Then the 
pure braid group may be realized as $P_n=G_n \rtimes \cdots \rtimes 
G_2$.  More generally, for $1\le\ll\le n$, the group 
$P_{n,\ll}=\ker(P_n\to P_\ll)$ may be realized as 
$P_{n,\ll}=G_n\rtimes \cdots\rtimes G_{\ll+1}$,\\ generated by 
$\c_{i,j}$, $\ll<j$.  Note that $P_n=P_{n,1}$.  For $\ll<j$, the 
monodromy homomorphisms $P_{j-1,\ll}\to\Aut(G_j)$ are given by the 
(restriction of the) Artin representation.  For $s<j$, 
we shall not distinguish 
between the braid $\c_{r,s}$ and the corresponding (right) 
automorphism $\c_{r,s}\in\Aut(G_j)$.  The action of $\c_{r,s}$ on 
$G_j$ is by conjugation: $\c_{r,s}(\c_{i,j}) = \c_{r,s}^{-1} \cdot 
\c_{i,j}^{} \cdot \c_{r,s}^{} = z_{i}^{} \cdot \c_{i,j}^{} \cdot 
z_{i}^{-1}$, where
\begin{equation} \label{eqn:Pnrels}
z_{i}=
\begin{cases}
\c_{r,j}\c_{s,j}&\text{if $i=r$ or $i=s$,}\\ \relax
[\c_{r,j},\c_{s,j}]&\text{if $r<i<s$,}\\
1&\text{otherwise.}
\end{cases}
\end{equation}
See 
Birman~\cite{Bi} 
and Hansen~\cite{Ha} for details, and as 
general references on braids.

\subsection{Some Fox Calculus}
\label{subsec:FoxCalc}
We first establish some notation and 
record some elements of the Fox Calculus \cite{Fox, Bi}, and 
results from \cite{CScc} necessary in the construction.  

Denote the integral group ring of a (multiplicative) group 
$G$ by $\Z{G}$.  We regard modules over $\Z{G}$ as 
left modules.  Elements of the free module $(\Z G)^{n}$ are viewed 
as row vectors, and $\Z G$-linear maps $(\Z G)^{n}\to (\Z G)^{m}$ 
are viewed as $n\times m$ matrices which act on the right.  For such 
a map $F$, denote the transpose by $F^\top$, and recall that $[F]^k$ 
denotes the map $\oplus_1^k F$.  Denote the $n\times n$ identity 
matrix by $\II_n$.  

For the single free group $G_j=\langle \c_{i,j}\rangle$, a free 
$\Z{G_j}$-resolution of $\Z$ 
is given by
\begin{equation} \label{eqn:oneres}
0 \to (\Z{G_j})^{j-1} \xrightarrow{\Delta_j} \Z{G_j} 
\xrightarrow{\epsilon} \Z
\to 0,
\end{equation}
where $\Delta_j=\begin{pmatrix}\c_{1,j}-1 &\cdots & \c_{j-1,j}-1 
\end{pmatrix}^{\top}$, and $\epsilon$ is the augmentation map, given 
by $\epsilon(\c_{i,j})=1$.  For each element $\c\in P_{j-1,\ll}$, 
conjugation by $\c$ induces an automorphism $\c:G_j\to G_j$, and a 
chain automorphism $\c_\*$ of \eqref{eqn:oneres}, which by the 
``fundamental formula of Fox Calculus,'' can be expressed as
\begin{equation} \label{eqn:ftc}
\begin{CD}
(\Z{G_j})^{j-1}        @>\Delta_j >>             \Z{G_j}\\
@VV\JJ(\c )\circ \tilde \c V   @VV\tilde \c V \\
(\Z{G_j})^{j-1}        @>\Delta_j >>             \Z{G_j}
\end{CD}
\end{equation}
where $\JJ(\c) = \left(\frac{\p\c(\c_{i,j})}{\p\c_{k,j}}\right)$ is 
the $(j-1)\times (j-1)$ Jacobian matrix of Fox derivatives of $\c$, 
and $\tilde\c$ denotes the extension of $\c$ to the group ring 
$\Z{G_j}$, resp.,~to $(\Z{G_j})^{j-1}$.  
For a second element $\b$ of $P_{j-1,\ll}$, we have
$(\c\cdot\b)_\* = (\b\circ\c)_\* = \b_\* \circ \c_\*$ by the 
``chain rule of Fox Calculus'': 
$\JJ(\b\circ\c)=\tilde\b(\JJ(\c))\cdot \JJ(\b)$.  In particular,
$\JJ(\c^{-1})=\tilde \c^{-1}(\JJ(\c)^{-1})$.

Now fix $\ll$, $1\le \ll \le n$, and consider the group 
$P_{n,\ll}=G_n\rtimes\cdots\rtimes G_{\ll+1}$.  Let 
$\RR=\Z{P_{n,\ll}}$ denote the integral group ring of $P_{n,\ll}$, For 
$\c\in P_{j-1,\ll}$ as above, define $m_\c:\RR\to \RR$ by 
$m_\c(r)=\c\cdot r$.  From \eqref{eqn:ftc} and extension of scalars, 
we obtain
\begin{equation*}
\begin{CD}
\RR\otimes _{\Z G_j} (\Z  G_j)^{j-1} @>{\id \otimes \Delta_j }>>
\RR\otimes _{\Z  G_j} \Z  G_j \\
@VVm_{\c} \otimes \JJ(\c )\circ \tilde \c V
@VVm_{\c} \otimes \tilde \c V \\
\RR\otimes _{\Z  G_j} (\Z  G_j)^{j-1} @>{\id \otimes \Delta_j }>>
\RR\otimes _{\Z  G_j} \Z  G_j \\
\end{CD}
\end{equation*}
The map $m_{\c} \otimes \JJ(\c )\circ \tilde \c $ and the canonical 
isomorphism $\RR\otimes _{\Z G_j} (\Z G_j)^{j-1} \cong \RR^{j-1}$ 
define an $\RR$-linear automorphism $\rho_j (\c ):\RR^{j-1} \to 
\RR^{j-1}$, whose matrix is $\c\cdot \JJ(\c)$, see 
\cite[Lemma~2.4]{CScc}.  Furthermore, we have the following.

\begin{lem}[{\cite[Lemma 2.6]{CScc}}]\label{lem:replem} 
For each $j$, $2\le j\le n$, the action of the group $P_{j-1,\ll}$ on 
the free group $G_j$ gives rise to a representation $\rho_{j}:
P_{j-1,\ll}\to \Aut _{\RR}(\RR^{j-1})$ with the property that 
$\rho _{j}(\c)=m_{\c}\otimes \JJ(\c)\circ \tilde \c$ for every $\c 
\in 
P_{j-1,\ll}$.
\end{lem}

\begin{rem} \label{rem:tilderem}
Via the convention $\rho _{j}(\c_{p,q}) = \II_{j-1}$ for $q\ge j$, 
the 
above extends to a representation $\rho _{j}:P_{n,\ll}\to \Aut 
_{\RR}(\RR^{j-1})$ of the entire group $P_{n,\ll}$.  We denote by   
$\tilde \rho _{j}:\RR\to \End _{\RR}(\RR^{j-1})$ the extension of 
$\rho _{j}$ to the group ring $\RR$.  We also use $\tilde \rho _{j}$ 
to denote the homomorphism $\homR _{\RR}(\RR^{m},\RR^{n})\to 
\homR_{\RR}(\RR^{m(j-1)},\RR^{n(j-1)})$ defined by replacing each 
entry $x$ of an $m \times n$ matrix by $\tilde \rho _{j}(x)$.
\end{rem}

\subsection{The Resolution}
\label{subsec:ZGres}
We now recall the construction of the free resolution $\epsilon 
:C_{\bullet }=C_{\bullet }(G)\to \Z $ over the ring $\RR=\Z G$ 
from~\cite{CScc}, in the case where $G=P_{n,\ll}$ is the fundamental 
group of the complement of the discriminantal arrangement 
$\A_{n,\ll}$.  If $J=\{j_1,\dots, j_q\}\subseteq [\ll+1,n]$, recall 
that for $p<q$, $J_p=\{j_1,\dots,j_p\}$ and $J^p=J\setminus J_p$.  For 
such a set $J$, let $C_q^J$ be a free $\RR$-module of rank 
$(j_1-1)\cdots(j_q-1)$.

Let $C_{0}=\RR$, and, for $1\le q\le n-\ll $, let $C_{q} = 
\bigoplus_{|J|=q } C_q^J$, where the sum is over all $J\subseteq 
[\ll+1,n]$.  The augmentation map, $\epsilon : C_{0}\to \Z $, is the 
usual augmentation of the group ring, given by $\epsilon (\c)=1$, for 
$\c\in P_{n,\ll}$.  We define the boundary maps of $C_{\bullet }$ by 
recursively specifying their restrictions $\Delta ^{J}$ to the 
summands $C_q^J$ as follows:
 
\smallskip

\noindent If $J=\{j\}$, we define $\Delta _J:C_1^{J}=\RR^{j-1}\to 
\RR=C_0$ as in the resolution \eqref{eqn:oneres}, by $\Delta _J =
\begin{pmatrix}\c_{1,j}-1&\cdots&\c_{j-1,j}-1 \end{pmatrix}^\top$. 

\smallskip

\noindent In general, if $J=\{j_1,\dots,j_q\}$, then 
$J^1=\{j_2,\dots,j_q\}$ and $J_{q-1}=\{j_1,\dots,j_{q-1}\}$, and we 
define $\Delta_J:C_q^J \to C_{q-1}^{J^1}$ by 
$\Delta_J=-\tilde\rho_{j_q}( \Delta_{J_{q-1}})$

\smallskip

\noindent Now define $\Delta ^{J}: C_q^J\to \bigoplus _{p=1}^{q} 
C_{q-1}^{J \setminus \{j_p\}}$ by 
\begin{equation*}
\Delta ^{J} =
\left (\Delta^{}_{J^{}},
\left [\Delta _{J^1}\right ]^{d_{1}},\dots ,
\left [\Delta _{J^{p}}\right ]^{d_p},
\dots , \left [\Delta _{J^{q-1}}\right ]^{d_{q-1}}\right ),
\end{equation*}
where $d_p=(j_1-1)\cdots (j_p-1)$.

\smallskip

\noindent Finally, define $\partial_q :C_{q}\to C_{q-1}$ by $\DS 
{\partial_{q} = \sum _{|J|=q } \Delta ^{J}}$.

\begin{thm}[{\cite[Theorem 2.10]{CScc}}]
\label{thm:resthm}
Let $\RR={\Z}P_{n,\ll}$ be the integral group ring of the group 
$P_{n,\ll}$.  Then the system of $\RR$-modules and homomorphisms 
$(C_\*,\partial_\*)$ is a finite, free resolution of $\Z $ over $\RR$.
\end{thm}

\begin{rem} \label{rem:mapcone1}
The proof of this result in \cite{CScc} makes use of a mapping cone 
decomposition of the complex $(C_\*,\partial_\*)$.  This 
decomposition 
may be described as follows.  Let $(\hC_\*,\hp_\*)$ denote the 
subcomplex of $(C_\*,\partial_\*)$ with terms $\hC_q = 
\bigoplus_{\ll+1\notin J} C_q^J$, and boundary maps 
$\hp_q=\partial_q\vert_{\hC_q}$ given by restriction.  The complex 
$\hC_\*$ may be realized as $\hC_\*=C_\*(P_{n,\ll+1}) 
\otimes_{P_{n,\ll}} \RR$, where $\epsilon:C_\*(P_{n,\ll+1})\to\Z$ is 
the resolution over $\Z{P_{n,\ll+1}}$ obtained by applying the above 
construction to the group $P_{n,\ll+1}< P_{n,\ll}$.

Let $(D_\*,\partial_\*^D)$ denote the direct sum of $\ll$ copies of 
the complex $\hC_\*$, with the sign of the boundary map reversed.  
That is, $D_q=(\hC_q)^\ll$ and $\partial_q^D=-[\hp_q]^\ll$.  The terms 
of this complex may be expressed as $D_q=\bigoplus_{\ll+1\in K} 
C_{q+1}^{K}$, where $|K|=q$.  Using this description, define a map 
$\Xi_\*:D_\*\to\hC_\*$ by setting the restriction of $\Xi_q$ to the 
summand $C_{q+1}^{K}$ of $D_q$ to be equal to $\Delta_K: C_{q+1}^{K} 
\to C_q^{J} \subset\hC_q$, where $K=\{\ll+1\}\cup J$.

As shown in \cite{CScc}, the map $\Xi_\*:D_\*\to\hC_\*$ is a chain 
map, and the original complex $(C_\*, \partial_\*)$ may be realized as 
the mapping cone of $\Xi_\*$.  Explicitly, the terms of $C_\*$ 
decompose as $C_q=D_{q-1}\oplus \hC_q$.  With respect to this 
decomposition, the boundary map $\partial_{q+1}: C_{q+1}\to C_q$ is 
given by $\partial_{q+1}(u,v)= 
(-\partial^D_q(u),\Xi_q(u)+\widehat\partial_{q+1}(v))$.
\end{rem}

\subsection{Rank One Local Systems}
\label{subsec:locsys}
The abelianization of the group $P_{n,\ll}$ is free abelian of rank 
$N=\binom{n}{2}-\binom{\ll}{2}$.  Let $\T$ denote the complex torus, 
with coordinates $t_{{i,j}}$, $\ll+1\le j\le n$, $1\le i<j$.  Each 
point $\bt \in \T$ gives rise to a rank one representation 
$\nu_\bt:P_{n,\ll}\to\C^*$, $\c_{{i,j}}\mapsto t_{{i,j}}$, an 
associated $P_{n,\ll}$-module $L=L_\bt$, and a rank one local system 
$\LL=\LL_\bt$ on the configuration space $M_{n,\ll}$.  The homology 
and cohomology of $P_{n,\ll}$ with coefficients in $L$ (resp.,~that of 
$M_{n,\ll}$ with coefficients in $\LL$) are isomorphic to the homology 
and cohomology of the complexes $\bC_\*:=C_\*\otimes_{P_{n,\ll}} L$ 
and $\bC^\*:=\Hom_{P_{n,\ll}}(C_\*,L)$ respectively, see \cite{Brown}.

The terms of these complexes, $\bC_q = C_q \otimes_{\RR} \C$ and 
$\bC^q=\text{Hom}_{P_{n,\ll}}(C_q,L)$, are finite dimensional complex 
vector spaces.  Notice that $\dim \bC_{q} = \dim \bC^{q} = \dim A^q = 
\sum_{|J|=q} (j_{1}-1)\cdots (j_{q}-1)$, where the sum is over all 
$J\subseteq [\ll+1,n]$.  Denote the boundary maps of $\bC_\*$ and 
$\bC^\*$ by $\bp_q(\bt):\bC_q\to\bC_{q-1}$ and 
$\bd^q(\bt):\bC^q\to\bC^{q+1}$ .  As we follow \cite{Brown} in our 
definition of $\bC^\*$, these maps are related by
\begin{equation} \label{eq:dual}
\bd^q(\bt)(u)(x) = (-1)^{q}u(\bp_{q+1}(\bt)(x)) 
\end{equation}
for $u\in \bC^q$ and $x\in C_{q+1}$.  To describe these maps further, 
we require some notation.

Consider the evaluation map $\RR\times \T\to \C$, which takes an 
element $f$ of the group ring, and a point $\bt$ in $\T$ and yields 
$f(\bt):=\tilde\nu_\bt(f)$, the evaluation of $f$ at $\bt$.  
Fixing $f\in \RR$ and allowing $\bt\in\T$ to vary, we get a 
holomorphic map $\mathbf{f}:\T\to\C$.  More generally, we have the map 
$\Mat_{r\times s}(\RR)\times \T\to \Mat_{r\times s}(\C)$, $(F,\bt) 
\mapsto F(\bt):=\tilde\nu_\bt(F)$.  For fixed $F\in\Mat_{p\times 
q}(\RR)$, we get a map $\mathbf{F}:\T\to\Mat_{r\times s}(\C)$.  With 
these conventions, if $\dim\bC_q=r$ and $\dim\bC_{q+1}=s$, the 
boundary maps of the complexes $\bC_{\*}$ and $\bC^\*$ may be viewed 
as evaluations, $\bp_q(\bt)$ and $\bd^q(\bt)$, of maps 
$\bp_q:\T\to\Mat_{r\times s}(\C)$ and $\bd^q:\T\to\Mat_{s\times 
r}(\C)$.

We shall subsequently be concerned with the derivatives of these maps 
at the identity element $\bone=(1,\dots,1)$ of $\T$.  The 
(holomorphic) tangent space of $H^1(M_{n,\ll};\C^*)=\T$ at $\bone$ is 
$H^1(M_{n,\ll};\C)=\C^N$, with coordinates $\l_{i,j}$.  The 
exponential map $\TT_\bone\T\to\T$ is the coefficient map 
$H^1(M_{n,\ll};\C) \to H^1(M_{n,\ll};\C^*)$ induced by 
$\exp:\C\to\C^*$, $\l_{i,j}\mapsto e^{\l_{i,j}}=t_{i,j}$.  For an 
element $f$ of $\RR$, the derivative of the corresponding map 
$\mathbf{f}:\T\to\C$ at $\bone$ is given by $\mathbf{f}_*:\C^N\to\C$, 
$\mathbf{f}_*(\bl) = {\frac{d}{dx}}\bigr\rvert _{x=0} f(\dots 
e^{x\l_{i,j}}\dots)$.  More generally, for $F\in \Mat_{r\times 
s}(\RR)$, we have $\mathbf{F}_*:\C^N \to \Mat_{r\times s}(\C)$.

\section{A Complex of Derivatives}
\label{sec:dercx}
We now relate the cohomology theories $H^*(A^\*,\mu^\*(\bl))$ and 
$H^*(M_n;\LL)$ by relating the complexes $(A^\*,\mu^\*(\bl))$ and 
$(\bC^\*,\bd^\*(\bt))$.  As above, let $(\bp_q)_*$ and $\bd^q_*$ 
denote the derivatives of the maps $\bp_q$ and $\bd^q$ at 
$\bone\in\T$.

\begin{thm} \label{thm:OS&CG}
The complex $(A^\*,\mu^\*(\bl))$ is a linear approximation of the 
complex $(\bC^\*,\bd^\*(\bt))$.  For each $\bl\in\C^N$, the system of 
complex vector spaces and linear maps $(\bC^\*,\bd^\*_*(\bl))$ is a 
complex.  For each $q$, we have $A^q\cong \bC^q$, and, under this 
identification, $\mu^q(\bl)=\bd^q_*(\bl)$.
\end{thm}

{F}rom the discussions in sections ~\ref{subsec:Hbraid} 
and~\ref{subsec:locsys}, it is clear that $A^q\cong \bC^q$.  In light 
of the sign conventions \eqref{eq:dual} used in the construction of 
the complex $(\bC^\*,\bd^\*(\bt))$ and the fact that 
$(A^\*,\mu^\*(\bl))$ is a complex, to show that 
$(\bC^\*,\bd^\*_*(\bl))$ is a complex, and to prove the theorem, it 
suffices to establish the following.

\begin{prop} \label{prop:maps}
For each $q$, we have 
$\mu^q(\bl)=(-1)^q\bigl[(\bp_{q+1})_*(\bl)\bigr]^\top$.
\end{prop}

The maps $\mu^q(\bl)$ were analyzed in 
section~\ref{subsec:wedge}.  We now carry out a similar analysis of 
the maps $(\bp_{q+1})_*(\bl)$.

\subsection{Some Calculus} \label{subsec:calc}

We first record some facts necessary for this analysis.  Recall that 
$\RR$ denotes the integral group ring of the group $P_{n,\ll}$.  For 
$f,g\in\RR$, the Product Rule yields $(\mathbf{f\cdot 
g})_*(\bl)=\mathbf{f}_*(\bl)\cdot \mathbf{g}(\bone)+ 
\mathbf{f}(\bone)\cdot \mathbf{g}_*(\bl)$.  Similarly, for $F \in 
\Mat_{p\times q}(\RR)$ and $G\in \Mat_{q\times r}(\RR)$, matrix 
multiplication and the differentiation rules yield 
\begin{equation}  
\label{eq:leibniz}
(\mathbf{F\cdot G})_*(\bl)=\mathbf{F}_*(\bl)\cdot\mathbf{G}(\bone)+
\mathbf{F}(\bone)\cdot\mathbf{G}_*(\bl).
\end{equation}
As an immediate consequence of the Product Rule, for $\c,\zeta\in 
P_{n,\ll}$ and $\tau=[\zeta,\c]$ a commutator, we have $(\bc^{-1})_* = 
-\bc_*$, and $\boldsymbol{\tau}_*=0$.  Consequently, $(\bzeta\cdot \bc 
\cdot \bzeta^{-1})_{*} = \bc_{*}$.

Now recall the representations $\rho_j$ defined in 
Lemma~\ref{lem:replem}, and used in the construction of the resolution 
$C_\*$.  Associated to each $\c \in P_{j-1,\ll}$, we have a map 
$\br_{j}(\c):\T \to \Aut(\C^{j-1})$.  Since $\c$ acts on the free 
group by conjugation, we have $\br_j(\c)(\bone)=\II_{j-1}$.  Identify 
$\End( \C^{j-1})$ as the tangent space to $\Aut(\C^{j-1})$ at the 
identity, and denote the derivative of the map $\br_{j}(\c)$ at 
$\bone$ by $\br_{j}(\c)_*:\C^N \to \End(\C^{j-1})$.

Define $(\br_{j})_{*}:P_{j-1,\ll}\to\text{Hom}(\C^{N},\End(\C^{j-1}))$ 
by $(\br_{j})_{*}(\c)=\br_{j}(\c)_*$.  The chain rule of Fox Calculus 
and a brief computation reveal that $(\br_{j})_{*}$ is a homomorphism, 
and is trivial on the commutator subgroup $P_{n,\ll}'$.  This yields a 
map $\C^{N}\to\text{Hom}(\C^{N},\End(\C^{j-1}))$, $\l_{r,s} \mapsto 
\br_j(\c_{r,s})_*$, which we continue to denote by $(\br_{j})_{*}$.  
For $\c\in P_{n,\ll}$, view the derivative, $\bc_*(\bl)=\sum c_{r,s} 
\l_{r,s}$, of the corresponding map $\bc$ as a linear form in the 
$\l_{r,s}$.  Then we have the following ``chain rule'':
\begin{equation} \label{eq:chainrule}
\br_j(\c)_*(\bl) = \sum c_{r,s} (\br_j)_*(\l_{r,s}) = 
(\br_j)_*(\bc_*(\bl)).
\end{equation}
In particular, $\br(\c_{r,s})_* = \br_*(\l_{r,s})$, which we now 
compute.
\begin{lem} \label{lem:onepurebraid}
For $r<s<j$, the derivative of the map $\br_j(\c_{r,s})$
is given by
\[
\br_j(\c_{r,s})_*(\bl)=
\begin{pmatrix}
\l_{r,s}\cdot\II_{r-1}&0&0&0&0\\
0        &\l_{r,s}+\l_{s,j}&0&-\l_{r,j}&0\\
0        &0      &\l_{r,s}\cdot\II_{s-r-1}&0&0\\
0        &-\l_{s,j}       &0        &\l_{r,s}+\l_{r,j}&0\\
0        &0           &0        &0  &\l_{r,s}\cdot\II_{j-s-1}
\end{pmatrix}\!.
\]
\end{lem}

\begin{proof} The matrix of $\rho_j(\c_{r,s})$ is $\c_{r,s}\cdot 
\JJ(\c_{r,s})$, where $\JJ(\c_{r,s})$ is the Fox Jacobian.  Thus, 
$\br_j(\c_{r,s})(\bt)= t_{r,s}\cdot \bJ(\c_{r,s})(\bt) = 
(t_{r,s}\cdot\II_{j-1})\cdot \bJ(\c_{r,s})(\bt)$, where 
$\bJ(\c_{r,s})(\bt)$ is the map induced by the Fox Jacobian.  By the 
Product Rule~\eqref{eq:leibniz}, we have
\[
\br_j(\c_{r,s})_*(\bl) = (\l_{r,s}\cdot\II_{j-1})\cdot
\bJ(\c_{r,s})(\bone) + \bJ(\c_{r,s})_*(\bl).
\]
The action of $\c_{r,s}$ on the free group $G_j=\langle 
\c_{i,j}\rangle$ is recorded in~\eqref{eqn:Pnrels}.  Computing Fox 
derivatives and evaluating at $\bt$ yields the familiar Gassner matrix 
of $\c_{r,s}$ (see~\cite{Bi}),
\[
\bJ(\c_{r,s})(\bt) = \begin{pmatrix}
\II_{r-1}&0&0&0&0\\
0        &1-t_{r,j}+t_{r,j}t_{s,j}&0&t_{r,j}(1-t_{r,j})&0\\
0        &\vec u      &\II_{s-r-1}&-\vec u&0\\
0        &1-t_{s,j}       &0        &t_{r,j}&0\\
0        &0           &0        &0  &\II_{j-s-1}
\end{pmatrix},
\]
where $\vec u = \begin{pmatrix} (1-t_{r+1,j})(1-t_{r,j})&\cdots& 
(1-t_{s-1,j})(1-t_{r,j}) \end{pmatrix}^\top$.  Since 
$\bJ(\c_{r,s})(\bone)=\II_{j-1}$, the result follows upon 
differentiating $\bJ(\c_{r,s})(\bt)$.
\end{proof}

\subsection{Boundary Map Derivatives} \label{subsec:Cmaps}
We now obtain an inductive formula for the derivatives of the boundary 
maps of the complex $(\bC_\*,\bp_\*(\bt))$.  The mapping cone 
decomposition of the resolution $(C_{\*},\partial_{\*})$ discussed in 
Remark~\ref{rem:mapcone1} gives rise to an analogous decomposition of 
the complex $(\bC_\*,\bp(\bt))$.  Specifically, the terms decompose as 
$\bC_{q}=\bD_{q-1}\oplus\widehat\bC_{q}$, and with respect to this 
decomposition, the matrix of the boundary map 
$\bp_{q+1}(\bt):\bC_{q+1}\to\bC_{q}$ is given by
\begin{equation} \label{eq:mapconedecomp}
\bp_{q+1}(\bt)=\begin{pmatrix} -\bp^D_q(\bt) & \bX_q(\bt)\\
0 & \widehat\bp_{q+1}(\bt) \end{pmatrix}.
\end{equation}
Up to sign, the complex $(\bD_{\*}, \bp^D_{\*}(\bt))$ is a direct sum 
of $\ll$ copies of the complex 
$(\widehat\bC_{\*},\widehat\bp_{\*}(\bt))$, which arises from the 
group $P_{n,\ll+1}<P_{n,\ll}$.  In light of this, we restrict our 
attention to the chain map $\Xi_{\*}$ and its components 
$\Delta_{\{\ll+1,J\}}$, their evaluations $\bX_\*(\bt)$ and 
$\bDD_{\{\ll+1,J\}}(\bt)$, and the derivatives of these evaluations at 
$\bone$.

For $J=\{j_1,\dots,j_q\}\subseteq [\ll+2,n]$, let 
$\rho^{}_{J}=\tr_{j_q} \circ \cdots \circ \tr_{j_1}$ and 
$d_J=(j_1-1)\cdots (j_q-1)$.  Then $\D_{\{\ll+1,J\}} = 
(-1)^q\rho^{}_{J} (\D_{\ll+1})$, where $\D_{\ll+1} = \begin{pmatrix} 
\c_{1,\ll+1}-1&\cdots&\c_{\ll,\ll+1}-1 \end{pmatrix}^\top$, and the 
matrix of $\D_{\{\ll+1,J\}}$ is $\ll\cdot d_J \times d_J$ with $d_J 
\times d_J$ blocks $(-1)^q\rho^{}_{J} (\c_{m,\ll+1}-1)$, $1\le m \le 
\ll$.  We concentrate our attention on one such block.

Fix $m$, $1\le m \le \ll$, and let $M$ denote the matrix of 
$\rho^{}_{J} (\c_{m,\ll+1}-1)$.  Similarly, let $M'$ denote the matrix 
of $\rho^{}_{J_{q-1}} (\c_{m,\ll+1}-1)$.  Then $M$ is the matrix of 
$\tr_{j_q}(M')$.  Since $M$ is $d_J\times d_J$, its rows and columns 
are naturally indexed by sets $R=\{r_1,\dots,r_q\}$ and 
$I=\{i_1,\dots,i_q\}$, $1\le r_p,i_p \le j_p-1$.  We thus denote the 
entries of $M$ by $M_{R,I}$.  With these conventions, we have 
\begin{equation} \label{eq:recursion}
M_{R,I} = \left[\tr_{j_q}(M'_{R',I'})\right]_{r_q,i_q}
\end{equation}
where, for instance, $I'=I_{q-1}=I\setminus\{i_q\}$.  

Now consider the block $\bM(\bt)$ of $\bDD_{\{\ll+1,J\}}(\bt)$ arising 
from the block $M$ of the matrix of $\D_{\{\ll+1,J\}}$ above.  Recall 
the notion of an $I$-admissible set from 
Definition~\ref{def:admissible}, and recall that $\epsilon_{R,I}=1$ if 
$R=I$, and $\epsilon_{R,I}=0$ otherwise.

\begin{thm} \label{thm:Centries}
Let $J \subseteq [\ll+2,n]$ and let $M$ denote the matrix of 
$\rho_J(\c_{m,\ll+1}-1)$.  Then the entries of the derivative, 
$\bM_*(\bl)$, of the evaluation $\bM(\bt)$ are given by
\[
[\bM_*(\bl)]_{R,I} = (-1)^{|R\setminus R\cap I|} \Bigl(
\epsilon_{R,I} \l_{m,\ll+1} + \sum_{j \in J} \sum_K \l_{k,j}\Bigr),
\]
where, if $j=j_p$, the second sum is over all $I_p$-admissible sets 
$K=\{k_{s_1},\dots,k_{s_t},k\}$ for which $R\setminus R\cap I 
\subseteq K$.
\end{thm}

\begin{proof}
The proof is by induction on $|J|$.

If $J=\{j\}$, then $M=\c_{m,\ll+1}\cdot \JJ(\c_{m,\ll+1}) - \II_{j-1}$ 
is the matrix of $\tr_j(\c_{m,\ll+1}-1)$, so $\bM(\bt) = 
\boldsymbol{\c}_{m,\ll+1}(\bt)\cdot \bJ(\c_{m,\ll+1})(\bt)-\II_{j-1}$.  
Since the derivative of the constant $\II_{j-1}$ is zero, the entries 
of $\bM_*(\bl)$ are given by Lemma~\ref{lem:onepurebraid} (with $r=m$ 
and $s=\ll+1$).  In this instance, we have $I=\{i\}$, and a set 
$K=\{k\}$ is $I$-admissible if $k\neq i$ and $\{k,i\}=\{m,\ll+1\}$.  
It follows that the case $|J|=1$ is a restatement of 
Lemma~\ref{lem:onepurebraid}.

In general, let $J=\{j_1,\dots,j_{q-1},j_q\}$ and write $M_{R,I} = 
\left[\tr_{j_q}(M'_{R',I'})\right]_{r_q,i_q}$ as in 
\eqref{eq:recursion} above.  Then we have
\[
[\bM(\bt)]_{R,I} = 
\left[\tilde{\boldsymbol{\rho}}_{j_q}(M'_{R',I'})(\bt)\right]_{r_q,i_q}
\text{\quad and\quad }
[\bM_*(\bl)]_{R,I} = 
\left[\tilde{\boldsymbol{\rho}}_{j_q}(M'_{R',I'})_*(\bl)\right]_{r_q,i_q}.
\]
By induction, the entries of the matrix $\bM'_*(\bl)$ are given by
\[
[\bM'_*(\bl)]_{R',I'} = (-1)^{|R'\setminus R'\cap I'|} \Bigl(
\epsilon_{R',I'} \l_{m,\ll+1} + \sum_{j_p \in J'} \sum_K 
\l_{k_p,j_p}\Bigr),
\]
where $J'=J_{q-1}=J\setminus\{j_q\}$, and for $j_p \in J'$, the second 
sum is over all $I_p$-admissible sets $K$ for which $R'\setminus 
R'\cap I' \subseteq K$.

By the chain rule~\eqref{eq:chainrule}, the entries of 
$\bM_*(\bl)$ are given by
\begin{equation} \label{eq:recursion2}
\begin{split}
\left[\bM_*(\bl)\right]_{R,I} &=  
\left[\tilde{\boldsymbol{\rho}}_{j_q}(M'_{R',I'})_*(\bl)\right]_{r_q,i_q} 
=
\left[(\boldsymbol{\rho}_{j_q})_*\bigl((\bM'_{R',I'})_*
(\bl)\bigr)\right]_{r_q,i_q}\\
&= 
\Bigl[S 
\Bigl(\epsilon_{R',I'}(\boldsymbol{\rho}_{j_q})_*(\l_{m,\ll+1})
+ \sum_{j_p \in J'} \sum_K 
(\boldsymbol{\rho}_{j_q})_*(\l_{k_p,j_p})
\Bigr)\Bigr]_{r_q,i_q},\\
\end{split}
\end{equation}
where $S=(-1)^{|R'\setminus R'\cap I'|}$.  
By Lemma~\ref{lem:onepurebraid}, for $r<s<j_q$, we have
\begin{equation} \label{eq:term1}
\left[(\boldsymbol{\rho}_{j_q})_*(\l_{r,s})\right]_{r_q,i_q} =
\begin{cases}
\l_{{r,s}}+ \l_{k_q,j_q} 
&\text{if $i_q=r_q$ and $\{k_q,i_q\}=\{r,s\}$,}\\
-\l_{r_q,j_q} &\text{if $i_q\neq r_q$ and $\{r_q,i_q\}=\{r,s\}$,}\\
0 &\text{otherwise.}
\end{cases}
\end{equation}

The entries of $\bM_*(\bl)$ may be calculated 
from~\eqref{eq:recursion2} using \eqref{eq:term1}, yielding the 
formula in the statement of the theorem.  We conclude the proof by 
making several observations which elucidate this calculation.

First consider the case $R'=I'$.  Then $S=1$ and $\epsilon_{R',I'}=1$.  
If $r_q=i_q$, then the first case of \eqref{eq:term1} yields a 
contribution of $\l_{{m,\ll+1}}+\l_{k_q,j_q}$ to 
$\left[\bM_*(\bl)\right]_{I,I}$, provided that 
$\{k_q,i_q\}=\{m,\ll+1\}$.  Note that this condition implies that the 
set $K=\{k_q\}$ is $I$-admissible (and that $k_q\neq i_q$).  Note also 
that in this instance we have $R=I$, $R\setminus R\cap I = \emptyset 
\subset K$, $\epsilon_{R,I}=1$, and $|R\setminus R\cap I|=0$.

If $R'=I'$ and $r_q \neq i_q$, then the second case of 
\eqref{eq:term1} contributes $-\l_{r_q,j_q}$ to 
$\left[\bM_*(\bl)\right]_{R,I}$ if $\{r_q,i_q\}=\{m,\ll+1\}$.  In this 
instance, the set $\{r_q\}$ is $I$-admissible, and since $R'=I'$ and 
$r_q \neq i_q$, we have $|R\setminus R\cap I|=1$.

For general $R'$ and $I'$, suppose that $S\cdot\l_{k_p,j_p}$ is a 
summand of $\left[ \bM'_*(\bl)\right]_{R',I'}$ for some $p\le q-1$.  
Then, by the inductive hypothesis, this summand arises from an 
$I_p$-admissible set $K=\{k_{s_1},\dots,k_{s_t},k_p\}$ with $R' 
\setminus R'\cap I' \subseteq K$.  If $r_q=i_q$, then the first case 
of \eqref{eq:term1} yields a contribution of 
$S\cdot(\l_{{k_p,j_p}}+\l_{k_q,j_q})$ to 
$\left[\bM_*(\bl)\right]_{R,I}$, provided that 
$\{k_q,i_q\}=\{k_p,j_p\}$.  For such $k_q$, it is readily checked that 
the set $K\cup\{k_q\}$ is $I$-admissible.  Also, since $r_q=i_q$, we 
have $R\setminus R\cap I=R'\setminus R'\cap I'\subseteq K$.

If, as above, $S\cdot\l_{k_p,j_p}$ is a summand of $\left[ 
\bM'_*(\bl)\right]_{R',I'}$ and $r_q\neq i_q$, then the second case of 
\eqref{eq:term1} contributes $-S\cdot\l_{r_q,j_q}$ to 
$\left[\bM_*(\bl)\right]_{R,I}$ provided $\{r_q,i_q\}=\{k_p.j_p\}$.  
In this instance, the set $K\cup\{r_q\}$ is $I$-admissible, and since 
$r_q \neq i_q$, we have $R\setminus R\cap I = (R'\setminus R'\cap 
I')\cup\{r_q\} \subseteq K\cup\{r_q\}$, and $|R\setminus R\cap 
I|=|R'\setminus R'\cap I'|+1$.

Applying these observations to \eqref{eq:recursion2} above 
completes the proof.
\end{proof}

\subsection{Proof of Proposition~\ref{prop:maps}}
We now use Theorems~\ref{thm:OSformula} and \ref{thm:Centries} to show 
that the differential of the complex $(A^\*,\mu^\*(\bl))$ is given by 
$\mu^q(\bl)=(-1)^q\bigl[(\bp_{q+1})_*(\bl)\bigr]^\top$, where 
$(\bp_q)_*(\bl)$ is the derivative of the boundary map of the complex 
$(\bC_\*,\bp_\*(\bt))$, thereby proving Proposition~\ref{prop:maps} 
and hence Theorem~\ref{thm:OS&CG} as well.

The proof is by induction on $d=n-\ll$, the cohomological dimension of 
the group $P_{n,\ll}$, (resp., the rank of the discriminantal 
arrangement $\A_{n,\ll}$).

In the case $d=1$, the complexes $A^\*$ and $\bC_\*$ are given by 
\[A^0\xrightarrow{\mu^0(\bl)} A^1\quad\text{and}\quad 
\bC_1\xrightarrow{\bp_1(\bt)} \bC_0\] 
respectively, where $A^0=\bC_0=\C$, $A^1=\oplus_{i<n}\C a_{i,n}$, and 
$\bC_1=\C^{n-1}$.  The boundary maps are $\mu^0(\bl):1\mapsto 
\sum_{i<n}\l_{i,n}\cdot a_{i,n}$ and $\bp_1(\bt) 
=\begin{pmatrix}t_{1,n}-1 &\cdots & t_{n-1,n}-1\end{pmatrix}^{\top}$.  
Identifying $A^1$ and $\bC_1$ in the obvious manner, we have 
$\mu^0(\bl)=(-1)^0\bigl[(\bp_{1})_*(\bl)\bigl]^\top$.

In the general case, we identify $A^q$ and $\bC_q$ in an analogous 
manner.  In particular, the rows and columns of the matrix of the 
boundary map $\bp_{q+1}(\bt):\bC_{q+1} \to \bC_q$ are indexed by basis 
elements $a_{I,J}$ of $A^{q+1}$ and $A^q$, or simply by the underlying 
sets $I$ and $J$, respectively.  To show that 
$\mu^q(\bl)=(-1)^q\bigl[(\bp_{q+1})_*(\bl)\bigr]^\top$, we make use of 
the decomposition of the complex $A^\*$ established in 
Proposition~\ref{prop:directsum}, and that of $\bC_\*$ stemming from 
the mapping cone decomposition of the resolution $C_\*$ described in  
Remark~\ref{rem:mapcone1}.  Recall from \eqref{eq:osmapmatrix} and 
\eqref{eq:mapconedecomp} that with respect to these decompositions, 
the boundary maps may be expressed as
\[
\mu^{q}(\bl) = \pmatrix \mu^{q}_{B}(\bl)&0\\
\Psi^q(\bl)&\widehat \mu^{q}(\bl)\endpmatrix
\quad\text{and}\quad
\bp_{q+1}(\bt) = \pmatrix -\bp_{q}^\bD(\bt)&\bX_q(\bt)\\
0&\widehat\bp_{q+1}(\bt)\endpmatrix.
\]

The maps $\widehat \mu^{q}(\bl)$ and $\widehat\bp_{q+1}(\bt)$ are the 
boundary maps of the complexes $\widehat A^\*$ and $\widehat\bC_\*$ 
arising from the cohomology algebra $A(\A_{n,\ll+1})$ and fundamental 
group $P_{n,\ll+1}$ of the complement of the discriminantal 
arrangement $\A_{n,\ll+1}$.  So by induction, we have $\widehat 
\mu^{q}(\bl)= (-1)^q\bigl[(\widehat\bp_{q+1})_*(\bl)\bigr]^\top$ for 
each $q$.  Since the complexes $B^\* \cong (\widehat A^\*)^\ll$ and 
$\bD_\* \cong (\widehat\bC_\*)^\ll$ decompose as direct sums, with 
boundary maps $\mu^q_B(\bl)=-[\widehat\mu^{q-1}(\bl)]^\ll$ and 
$\bp^\bD_q(\bt)=-[\widehat\bp_q(\bt)]^\ll$ the inductive hypothesis 
also implies that
\[
\mu^q_B(\bl)= 
-\Bigl[(-1)^{q-1}\bigl[(\widehat\bp_q)_*(\bl)\bigr]^\top\Bigr]^\ll = 
(-1)^{q}\Bigl[\bigl[(\widehat\bp_q)_*(\bl)\bigr]^\ll\Bigr]^\top = 
(-1)^q\bigl[-(\bp^\bD_q)_*(\bl)\bigr]^\top.
\]

Thus it remains to show that 
$\Psi^q(\bl)=(-1)^q\bigl[(\bX_q)_*(\bl)\bigr]^\top$.  For this, it 
suffices to show that the restriction $\Psi^q_J(\bl):A^q_J \to 
A^{q+1}_{\{\ll+1,J\}}$ of $\Psi^q(\bl)$ is dual to the derivative of 
the summand $\bDD_{\{\ll+1,J\}}(\bt):\bC^{\{\ll+1, J\}}_{q+1} \to 
\bC^J_q$ of $\bX_q(\bt)$ for each $J=\{j_1,\dots,j_q\}\subseteq 
[\ll+2,n]$.  As noted in \eqref{eq:psiblocks}, the matrix of 
$\Psi^q_J(\bl)$ is $d_J \times \ll\cdot d_J$ with $d_J\times d_J$ 
blocks $\pi_{m,\ll+1}\circ \Psi^q_J(\bl)$, where 
$d_J=(j_1-1)\cdots(j_q-1)$.  Similarly, from the discussion in section 
\ref{subsec:Cmaps}, we have that the matrix of 
$\bDD_{\{\ll+1,J\}}(\bt)$ is $\ll\cdot d_J \times d_J$ with $d_J\times 
d_J$ blocks $(-1)^q\bigl(\br_J(\c_{m,\ll+1})(\bt) - \II_{d_J}\bigr)$.  
Comparing the formulas obtained in Theorem~\ref{thm:OSformula} and 
Theorem~\ref{thm:Centries}, we see that $\pi_{m,\ll+1}\circ 
\Psi^q_J(\bl) = \bigl[\br_J(\c_{m,\ll+1})_*(\bl)\bigr]^\top$.  It 
follows readily that $\Psi^q_J(\bl) = 
(-1)^q\bigl[\bigr(\bDD_{\{\ll+1,J\}}\bigr)_*(\bl)\bigr]^\top$, 
completing the proof.

\section{Cohomology Support Loci and Resonant Varieties}
\label{sec:loci}
In an immediate application of Theorem~\ref{thm:OS&CG}, we establish 
the relationship between the cohomology support loci of the complement 
of the discriminantal arrangement $\A_{n,\ll}$ and the resonant 
varieties of its Orlik-Solomon algebra.

Recall that each point $\bt\in\T$ gives rise to a local system 
$\LL=\LL_{\bt}$ on the complement $M_{n,\ll}$ of the arrangement 
$\A_{n,\ll}$.  For sufficiently generic $\bt$, the cohomology 
$H^{k}(M_{n,\ll},\LL_{\bt})$ vanishes (for $k<n-\ll$), see for 
instance \cite{Ko,CScc}.  Those $\bt$ for which 
$H^{k}(M_{n,\ll};\LL_\bt)$ does not vanish comprise the cohomology 
support loci 
\[
\Sigma^{k}_m(M_{n,\ll})=\{\bt\in\T \mid \dim H^{k}(M_{n,\ll}; 
\LL_{\bt})\ge m\}.
\]
These loci are algebraic subvarieties of $\T$, which, since 
$M_{n,\ll}$ is a $K(P_{n,\ll},1)$-space, are invariants of the group 
$P_{n,\ll}$.

Similarly, each point $\bl\in\C^{N}$ gives rise to an element 
$\omega=\omega_{\bl}\in A^{1}$ of the Orlik-Solomon algebra of the 
arrangement $\A_{n,\ll}$.  For sufficiently generic $\bl$, the 
cohomology $H^{k}(A^{\*},\mu^{\*}(\bl))$ vanishes (for $k<n-\ll$), see 
\cite{Yuz,Fa}.  Those $\bl$ for which $H^{k}(A^{\*},\mu^{\*}(\bl))$ 
does not vanish comprise the resonant varieties
\[
\RR^{m}_k(A)=\{\bl\in\C^{N} \mid \dim H^{k}(A^{\*},\mu^{\*}(\bl))\ge 
m\}.
\]
These subvarieties of $\C^{N}$ are invariants of the Orlik-Solomon 
algebra $A$ of $\A_{n,\ll}$.

Recall that $\bone=(1,\dots,1)$ denotes the identity element of $\T$.

\begin{thm} \label{thm:tcone}
Let $\A_{n,\ll}$ be a discriminantal arrangement with complement 
$M_{n,\ll}$ and Orlik-Solomon algebra $A$.  Then for each $k$ and each 
$m$, the resonant variety $\RR_{k}^{m}(A)$ coincides with the tangent 
cone of the cohomology support locus $\Sigma_{m}^{k}(M_{n,\ll})$ at 
the point $\bone$.
\end{thm}
\begin{proof}
For each $\bt\in\T$, the cohomology of $M_{n,\ll}$ with coefficients 
in the local system $\LL_{\bt}$ is isomorphic to that of the complex 
$(\bC^{\*},\bd^{\*}(\bt))$.  So $\bt \in \Sigma_{m}^{k}(M_{n,\ll})$ if 
and only if $\dim H^{k}(\bC^{\*},\bd^{\*}(\bt) \ge m$.  An exercise in 
linear algebra shows that
\[
\Sigma_{m}^{k}(M_{n,\ll})=\{\bt\in\T\mid
\rank \bd^{k-1}(\bt) + \rank \bd^{k}(\bt) \le \dim \bC^k-m\}.
\]

For $\bl\in \C^{N}$, we have $\bl\in \RR_{k}^{m}$ if $\dim 
H^{k}(A^{\*},\mu^{\*}(\bl)) \ge m$.  So, as above, 
\[\RR_{k}^{m}(A)=\{\bl\in\C^{N}\mid \rank \mu^{k-1}(\bl) + \rank 
\mu^{k}(\bl) \le \dim A^{k}-m\}.
\]
By Theorem~\ref{thm:OS&CG}, $\dim A^{k}=\dim \bC^{k}$ and 
$\mu^{k}(\bl)=\bd^{k}_{*}(\bl)$ for each $k$.  
Thus, 
\[\RR_{k}^{m}(A)=\{\bl\in\C^{N}\mid \rank \bd^{k-1}_{*}(\bl) + 
\rank \bd^{k}_{*}(\bl) \le \dim \bC^k-m\},
\]
and the result follows.
\end{proof}

The cohomology support loci are known to be unions of 
torsion-translated subtori of $\T$, see \cite{Ara}.  In particular, 
all irreducible components of $\Sigma_{m}^{k}(M_{n,\ll})$ passing 
through $\bone$ are subtori of $\T$.  Consequently, all irreducible 
components of the tangent cone are linear subspaces of $\C^{N}$.  So 
we have the following.

\begin{cor} \label{cor:falkconj}
For each $k$ and each $m$, the resonant variety $\RR^{m}_{k}(A)$ is 
the union of an arrangement of subspaces in $\C^{N}$.
\end{cor}

\begin{rem} For $k=1$, Theorem \ref{thm:tcone} and Corollary 
\ref{cor:falkconj} hold for any arrangement $\A$, see 
\cite{CScv,L3,LY}.  In particular, as conjectured by Falk 
\cite[Conjecture~4.7]{Fa}, the resonant variety $\RR_{1}^{m}(A(\A))$ 
is the union of a subspace arrangement.  Thus, Corollary 
\ref{cor:falkconj} above may be viewed as resolving positively a 
strong form of this conjecture in the case where $\A=\A_{n,\ll}$ is a 
discriminantal arrangement.
\end{rem}

\bibliographystyle{amsalpha}

\end{document}